\documentclass[12pt]{amsart}
\usepackage{amsmath,amssymb,amsfonts,hyperref}
\usepackage[margin=1in]{geometry}

\usepackage[T1]{fontenc}
\usepackage{color}
\newtheorem{theorem}{Theorem}

\newtheorem{corollary}{Corollary}

\newtheorem{remark}{Remark}

\theoremstyle{definition}

\newtheorem{example}{Example}

\def\<{\langle}
\def\>{\rangle}

\DeclareMathOperator{\cpct}{Cap}

\DeclareMathOperator{\dom}{dom}

\begin{document}
\title[Probabilistic characterizations of essential self-adjointness]{Probabilistic characterizations of essential self-adjointness and removability of singularities}
\author[Michael Hinz]{Michael Hinz$^1$}
\address{$^1$ Department of Mathematics,
University of Bielefeld, 33501 Bielefeld, Germany}
\email{mhinz@math.uni-bielefeld.de}
\author[Seunghyun Kang]{Seunghyun Kang$^2$}
\address{$^2$ Department of Mathematical Sciences, Seoul National University, GwanAkRo 1, Gwanak-Gu, Seoul 08826, Korea}
\email{power4454@snu.ac.kr}
\author[Jun Masamune]{Jun Masamune$^3$}
\address{$3$Department of Mathematics, Hokkaido University, Kita 10, Nishi 8, Kita-Ku, Sapporo, Hokkaido, 060-0810, Japan}
\email{jmasamune@math.sci.hokudai.ac.jp}

\begin{abstract}
We consider the Laplacian and its fractional powers of order less than one on the complement $\mathbb{R}^d\setminus\Sigma$ of a given compact set $\Sigma\subset \mathbb{R}^d$ of zero Lebesgue measure. Depending on the size of $\Sigma$, the operator under consideration, equipped with the smooth compactly supported functions on $\mathbb{R}^d \setminus \Sigma$, may or may not be essentially self-ajoint. We survey well known descriptions for the critical size of $\Sigma$ in terms of capacities and Hausdorff measures. In addition, we collect some known results for certain two-parameter stochastic processes. What we finally want to point out is, that, although a priori essential self-adjointness is not a notion directly related to classical probability, it admits a characterization via Kakutani type theorems for such processes.  
\tableofcontents

\end{abstract}


\maketitle


\section{Introduction}

In this note we would like to point out an interesting connection between some traditional and well-studied notions in analysis and an interesting, but perhaps slightly less known area in probability theory. More precisely, we outline the relation between uniqueness questions for self-adjoint extensions of the Laplacian and its powers on the one hand and hitting probabilities for certain two-parameter stochastic processes on the other. Although both, the analytic part and the probabilistic part of the results stated below are well-established, it seems that the existing literature did never merge these two different aspects. 

Recall that if a symmetric operator in a Hilbert space, considered together with a given dense initial domain, has a unique self-adjoint extension, then it is called \emph{essentially self-adjoint}. The question of essential self-adjointness has strong physical relevance, because the evolution of a quantum system is described in terms of a unitary group, the generator of a unitary group is necessarily self-adjoint, and different self-adjoint operators determine different unitary groups, i.e. different physical dynamics. See for instance \cite[Section X.1]{RS80}. Self-adjointness, and therefore also essential self-adjointness, are notions originating from quantum mechanics.

A related notion of uniqueness comes up in probability theory, more precisely, in the theory of Markov semigroups. Recall that any non-positive definite self-adjoint operator $L$ on a Hilbert space $H$ is uniquely associated with a non-negative definite closed and densely defined symmetric bilinear form $Q$ on $H$ by $Q(u,v)=-\left\langle u,Lv\right\rangle_H$, \cite[Section VIII.6]{RS80}, where $\left\langle\cdot,\cdot\right\rangle_H$ denotes the scalar product in $H$ and $u$ and $v$ are arbitrary elements of the domain of $Q$ and the domain of $L$, respectively. Now assume that $H$ is an $L^2$-space of real-valued (classes of) functions. Then, if for any $u$ from the domain of $Q$ also $|u|$ is in the domain of $Q$ and we have $Q(|u|,|u|)\leq Q(u,u)$, the form $Q$ is said to satisfy the \emph{Markov property}. In this case it is called a \emph{Dirichlet form}, and $L$ is the infinitesimal generator of a uniquely determined strongly continuous semigroup of symmetric Markov operators on $H$, sometimes also called a \emph{Markov generator}, \cite{BH91, Da89, FOT94}. We say that a non-positive definite symmetric operator in an $L^2$-space $H$, together with a given dense initial domain, is \emph{Markov unique}, if it has a unique self-adjoint extension in $H$ that generates a Markov semigroup. Different Markov generators determine different Markov semigroups and (disregarding for a moment important issues of  construction and regularity) this means that they define different Markov processes. So the notion of Markov uniqueness belongs to probability theory. It has strong relevance in the context of
classical mechanics and statistical physics. 

For a non-positive definite densely defined symmetric operator on an $L^2$-space essential self-adjointness implies Markov uniqueness, but the converse implication is false, see Examples \ref{Ex:exMu} and \ref{Ex:exesa} below or \cite{Takeda92}. Even if an operator is Markov unique, it may still have other self-adjoint extensions that do not generate Markov semigroups. It is certainly fair to say that a priori the notion of essential self-adjointness is a not a probabilistic notion. However, and this is what we would like to point out here, in certain situations essential self-adjointness can still be characterized in terms of classical probability. 

We consider specific exterior boundary value problems in $\mathbb{R}^d$. It is well-known that the Laplacian $\Delta$, endowed with the initial domain $C_0^\infty(\mathbb{R}^d)$ of smooth compactly supported functions on $\mathbb{R}^d$, has a unique self-adjoint extension in $L^2(\mathbb{R}^d)$. This unique self-adjoint extension is given by $(\Delta, H^2(\mathbb{R}^d))$, where given $\alpha>0$, the symbol $H^\alpha(\mathbb{R}^d)$ denotes the Bessel potential space of order $\alpha$, see Section \ref{S:Prelim} below. Similarly, the fractional Laplacians $-(-\Delta)^{\alpha/2}$ of order $\alpha>0$, endowed with the domain $C_0^\infty(\mathbb{R}^d)$, have unique self-adjoint extensions, respectively, namely $(-(-\Delta)^{\alpha/2}, H^{\alpha}(\mathbb{R}^d))$. In the present note we focus on the cases $0<\alpha\leq 2$. 

Given a compact set $\Sigma\subset \mathbb{R}^d$ of zero $d$-dimensional Lebesgue measure, we denote its complement by $N:=\mathbb{R}^d\setminus\Sigma$.
For any $0<\alpha\leq 2$ the operator $(-(-\Delta)^{\alpha/2}, C_0^\infty(N))$ is non-positive definite and symmetric on $L^2(N)=L^2(\mathbb{R}^d)$. We are interested in conditions on the size of $\Sigma$ so that $(-(-\Delta)^{\alpha/2}, C_0^\infty(N))$ is essentially self-adjoint. Of course 
one possible self-adjoint extension is the global operator $(-(-\Delta)^{\alpha/2}, H^\alpha(\mathbb{R}^d))$, which 'ignores' $\Sigma$. If $\Sigma$ is 'sufficiently small', it will not be seen, and there is no other self-adjoint extension. If $\Sigma$ is 'too big', it will registered as a boundary, leading to a self-adjoint extension different from the global one. 

As mentioned, the analytic background of this problem is classical and can for instance be found in the textbooks \cite{AH96, FOT94, Mazja}. See in particular \cite[Sections 13.3 and 13.4]{Mazja}.
For integer powers of the Laplacian on $\mathbb{R}^d$ a description of the critical size of $\Sigma$ in terms of capacities and Hausdorff measures had been given in \cite[Section 10]{AGMST13}, and to our knowledge this was the first reference that gave such a characterization of essential self-adjointness. For fractional powers such descriptions do not seem to exist in written form. A probabilistic description for the critical size of $\Sigma$, which we could not find anywhere in the existing literature, can be given in terms of suitable two-parameter processes as for instance studied in \cite{Hirsch95, Khoshnevisan02, KhoshnevisanXiao04, KhoshnevisanXiao05}. In essence, these descriptions are straightforward applications of Kakutani type theorems for multiparameter processes, see for instance \cite[Chapter 11, Theorems 3.1.1 and 4.1.1]{Khoshnevisan02}. In fact, using processes with more than two parameters one could even extend this type of results to fractional Laplacians of arbitrary order.
A philosophically related idea, namely a connection between Riesz capacities and the hitting behaviour of certain one-parameter Gaussian processes (that are not Markov processes except in the Brownian case) had already been studied in \cite{Kahane83}.

We would like to announce related forthcoming results for Laplacians on complete Riemannian manifolds, \cite{HinzMasamune}. An analytic description of essential self-adjointness for the Laplacian via capacities reads as in the Euclidean case, instead of traditional arguments for Euclidean spaces based on convolutions, \cite{AH96}, our proof uses the regularity theory for the Laplacian on manifolds, \cite{Grigoryan2009}, and basic estimates on the gradients of resolvent densities, \cite[Section 4.2]{Aubin82}. To proceed to a geometric description we use asymptotics of the resolvent densities, they are basically the same as those for Green functions, see for instance \cite[Section 4.2]{Aubin82}, \cite[Section 4.2]{Grigoryan1999} or \cite[Section 4.2]{LT87}. For a probabilistic description we restrict ourselves, at least for the time being, to the case of Lie groups. In this case we can still work with relatively simple two-parameter processes and use the potential developed in \cite{Hirsch95, HirschSong95} to connect them to capacities and essential self-adjointness. In the case of general complete Riemannian manifolds
one first has to raise the quite non-trivial question what could be suitable two-parameter processes taking values in manifolds, see the comments in \cite{HinzMasamune}.

A subsequent idea to be addressed in the near future concerns details of the relationship between stochastic processes and specific boundary value problems. For many interesting cases it is well understood how boundary value problems (such as Dirichlet, Neumann or mixed), encoded in the choice of domain for the associated Dirichlet form, determine the behaviour of associated one-parameter Markov processes. It would be interesting to see whether, and if yes, in what sense, the behaviour of related two-parameter processes can reflect given boundary value problems for the Laplacian, encoded in the choice of its domain as a self-adjoint operator.

In the next section we collect some preliminaries. In Section \ref{S:Muesa} we discuss analytic characterizations of Markov uniqueness and essential self-adjointness for fractional Laplacians. In Section \ref{S:geo} we provide geometric descriptions, and in Section \ref{S:Probab} we give probabilistic characterizations in terms of hitting probabilities for two-parameter processes.

\subsection*{Acknowledgements} We would like to thank Professors Masayoshi Takeda, Sergio Albeverio and Hiroaki Aikawa for helpful and inspiring discussions on the subject.

\section{Bessel potential spaces, capacities and kernels}\label{S:Prelim}

We provide some preliminaries on function spaces, fractional Laplacians, related capacities and kernels. Our exposition mainly follows \cite[Chapters 1-3]{AH96}. Given $\alpha>0$ we define the \emph{Bessel potential space of order $\alpha$} by
\[H^\alpha(\mathbb{R}^d)=\left\lbrace u \in L^2(\mathbb{R}^d): (1+|\xi|^2)^{\alpha/2}\hat{u}\in L^2(\mathbb{R}^d)\right\rbrace, \] 
where $u\mapsto\hat{u}$ denotes the Fourier transform of $u$. Together with the norm 
\[\left\|u\right\|_{H^\alpha(\mathbb{R}^d)}=\left\|(1+|\xi|^2)^{\alpha/2}\hat{u}\right\|_{L^2(\mathbb{R}^d)}\]
it becomes a Hilbert space. See for instance \cite{AH96, Mazja, Tri78, Tri97}. Using the fact that 
\[-\Delta f= (|\xi|^2 \hat{f})^\vee\]
for any $f\in\mathcal{S}(\mathbb{R}^d)$, where $\mathcal{S}(\mathbb{R}^d)$ denotes the space of Schwartz functions on $\mathbb{R}^d$ and $u\mapsto\check{u}$ the inverse Fourier transform, we can easily see that $(\Delta, C_0^\infty(\mathbb{R}^d))$ is essentially self-adjoint on $L^2(\mathbb{R}^d)$ with the unique self-adjoint extension $(\Delta, H^2(\mathbb{R}^d))$, see  for instance \cite[Theorem 3.5.3]{Da95}. For $\alpha>0$ we can define the \emph{fractional Laplacians $-(-\Delta)^{\alpha/2}$ of order $\alpha/2$} in terms of Fourier transforms by \[(-\Delta)^{\alpha/2}f=(|\xi|^{\alpha} \hat{f})^\vee.\]
Again it is not difficult to show that $(-\Delta)^{\alpha/2}$, endowed with the domain $C_0^\infty(\mathbb{R}^d)$, has a unique self-adjoint extension, namely $((-\Delta)^{\alpha/2}, H^{\alpha}(\mathbb{R}^d))$. One can proceed similarly as in \cite[Theorem 3.5.3]{Da95}, see also \cite[Theorem 1.2.7 and Lemma 1.3.1]{Da95}.

Given $\alpha>0$, we write 
\begin{equation}\label{E:Besselkernel}
\gamma_\alpha:=((1+|\xi|^2)^{-\alpha/2})^\vee
\end{equation}
to denote the \emph{Bessel kernel of order $\alpha$} and $\mathcal{G}_\alpha f:=\gamma_\alpha\ast f$ to denote the \emph{Bessel potential operator $\mathcal{G}_\alpha$ of order $\alpha$}, which defines a bijection from $\mathcal{S}(\mathbb{R}^d)$ into itself and also a bounded linear operator $\mathcal{G}_\alpha: L^2(\mathbb{R}^d)\to L^2(\mathbb{R}^d)$. In both interpretations we have $\mathcal{G}_\alpha=(I-\Delta)^{-\alpha/2}$. The image $\mathcal{G}_\alpha f$ of a measurable function $f:\mathbb{R}^d\to[0,+\infty]$ is a lower semicontinuous nonnegative function on $\mathbb{R}^d$, see \cite[Proposition 2.3.2]{AH96}. This implies that for any $f\in L^2_+(\mathbb{R}^d)$, where the latter symbol denotes the cone of nonnegative elements in $L^2(\mathbb{R}^d)$, its image $\mathcal{G}_\alpha f$ is a $[0,+\infty]$-valued function on $\mathbb{R}^d$, i.e. defined for any $x\in\mathbb{R}^d$. We can therefore define the \emph{$\alpha,2$-capacity $\cpct_{\alpha,2}(E)$} of a set $E\subset \mathbb{R}^d$ by
\[\cpct_{\alpha,2}(E)=\inf\left\lbrace \left\|f\right\|_{L^2(\mathbb{R}^d)}^2: \text{$f\in L^2_+(\mathbb{R}^d)$ and $\mathcal{G}_\alpha f(x)\geq 1$ for all $x\in E$}\right\rbrace,\]
with the convention that $\cpct_{\alpha,2}(E)=+\infty$ if no such $f$ exists, see \cite[Definition 2.3.3]{AH96}.

There is another, 'more algebraic' definition of a $\alpha,2$-capacity. For a compact set $K\subset\mathbb{R}^d$, define
\begin{multline}\label{E:algcap}
\cpct_{\alpha,2}'(K)=\inf\left\lbrace \left\|\varphi\right\|_{H^\alpha(\mathbb{R}^d)}^2: \text{ $\varphi\in C_0^\infty(\mathbb{R}^d)$ such that $\varphi(x)=1$}\right.\\
\left.\text{ for all $x$ from a neighborhood of $K$}\right\rbrace.
\end{multline}
Exhausting open sets by compact ones and approximating arbitrary sets from outside by open ones, this definition can be extended in a consistent manner to arbitrary subsets of $\mathbb{R}^d$. Now it is known that there exist constants $c_1,c_2>0$ such that for any compact set $K\subset \mathbb{R}^d$, we have
\begin{equation}\label{E:capscoincide}
c_1\:\cpct_{\alpha,2}(K)\leq \cpct_{\alpha,2}'(K)\leq c_2\:\cpct_{\alpha,2}(K),
\end{equation}
see \cite[Theorem 3.3]{Mazja72} for integer $\alpha$ and \cite[Section 2.7 and Corollary 3.3.4]{AH96} or \cite[Theorem A]{AP73} for general $\alpha$. We would like to remark that (\ref{E:capscoincide}) is based on certain truncation results for potentials. For $0<\alpha\leq 1$ 
the spaces $H^\alpha(\mathbb{R}^d)$ are domains of Dirichlet forms so that truncation properties 
are immediate from the Markov property. However, for $\alpha>1$ one needs to invest additional arguments, see for instance \cite[Sections 3.3, 3.5 and 3.7]{AH96}. 

As before, let $\alpha>0$. We say that a Radon measure $\mu$ on $\mathbb{R}^d$ has \emph{finite $\alpha$-energy} if 
\[\int_{\mathbb{R}^d} |v|d\mu\leq c\left\|v\right\|_{H^\alpha(\mathbb{R}^d)}\quad\text{ for all $v\in C_0^\infty(\mathbb{R})$.}\]
For a measure $\mu$ having finite $\alpha$-energy we can find a function $U^\alpha\mu\in H^\alpha(\mathbb{R}^d)$ such that 
\begin{equation}\label{E:finiteenergyint}
\left\langle U^\alpha\mu, v\right\rangle_{H^\alpha(\mathbb{R}^d)}=\int_{\mathbb{R}^d} v\:d\mu\quad\text{ for all $v\in \mathcal{S}(\mathbb{R}^d)$.}
\end{equation}
Using Fourier transforms this seen to be equivalent to requiring
\[\left\langle (1+|\xi|^2)^{\alpha} \widehat{U^\alpha\mu}, \hat{v}\right\rangle_{L^2(\mathbb{R}^d)}=\hat{\mu}(\widehat{v(-\cdot)})\quad\text{ for all $v\in \mathcal{S}(\mathbb{R}^d)$,}\]
what implies that $\widehat{U^\alpha \mu}=(1+|\xi|^2)^{-\alpha}\hat{\mu}$ in the sense of Schwartz distributions, and finally,
\[U^\alpha\mu=\gamma_{2\alpha}\ast\mu.\]
Note that by (\ref{E:Besselkernel}) we have
\begin{equation}\label{E:convo}
\gamma_{2\alpha}=\gamma_{\alpha}\ast\gamma_{\alpha}.
\end{equation}
We can define the $\alpha$-energy of $\mu$ as
\[E_\alpha(\mu):=\int_{\mathbb{R}^d}U^\alpha\mu\:d\mu,\]
and by (\ref{E:finiteenergyint}) this can be seen to equal $\left\|U^\alpha\mu\right\|_{H^\alpha(\mathbb{R}^d)}^2$. There is a dual definition of the $\alpha,2$-capacity: For a compact set $K\subset\mathbb{R}^d$ we have
\begin{equation}\label{E:dualdef}
\cpct_{\alpha,2}(K)=\sup\left\lbrace \frac{\mu(K)^2}{E_\alpha(\mu)}:\text{ $\mu$ is a Radon measure on $K$}\right\rbrace
\end{equation}
with the interpretation $\frac{1}{\infty}:=0$, see \cite[Theorem 2.2.7]{AH96}.

We finally collect some well-known asymptotics of the Bessel kernels. For $0<\alpha<d$ we have
\begin{equation}\label{E:asympzero}
\gamma_\alpha\sim c_{d,\alpha}|x|^{\alpha-d}\quad\text{ as $|x|\to 0$}
\end{equation}
with a positive constant $c_{d,\alpha}$ depending only on $d$ and $\alpha$, and for the limit case $\alpha=d$,
\begin{equation}\label{E:asympzerolog}
\gamma_d(x)\sim c_d (-\log|x|)\quad\text{ as $|x|\to 0$}
\end{equation}
with a positive constant $c_d$ depending on only on $d$. Moreover, it is known that
\begin{equation}\label{E:asympinfty}
\gamma_\alpha(x)=O(e^{-c|x|}) \quad\text{ as $|x|\to \infty$}.
\end{equation}
By (\ref{E:Besselkernel}) we have $\hat{\gamma}_\alpha(\xi)\leq |\xi|^{-\alpha}$ for all sufficiently large $\xi\in\mathbb{R}^d$. In the case $d<\alpha$ we therefore see that
the Bessel kernel $\gamma_\alpha$ is an element of $L^1(\mathbb{R}^d)$ and equals 
\begin{equation}\label{E:intcase}
\gamma_\alpha(x)=\int_{\mathbb{R}^d}\frac{e^{i\left\langle x,\xi\right\rangle}}{(1+|\xi|^2)^{\alpha/2}}\:d\xi,\quad x\in\mathbb{R}^d.
\end{equation}
See \cite[Sections 1.2.4 and 1.2.5]{AH96}.

\section{Markov uniqueness, essential self-adjointness and capacities}\label{S:Muesa}

Recall that $\Sigma\subset \mathbb{R}^d$ is a given compact set of zero Lebesgue measure and $N:=\mathbb{R}^d\setminus\Sigma$. We first state a well-known known result on Markov uniqueness. 
Using the definition (\ref{E:algcap}) of capacities together with traditional approximation arguments, which we will formulate below for the question of essential self-adjointness, one can obtain the following.

\begin{theorem}\label{T:Mu} Let $0<\alpha\leq 2$.
The fractional Laplacian $((-\Delta)^{\alpha/2}, C_0^\infty(N))$ is Markov-unique if and only if $\cpct_{\alpha/2,2}(\Sigma)=0$.
\end{theorem}

A classical guiding example for the case $\alpha=2$ is the following, which will be complemented for the cases $0<\alpha<2$ in Section \ref{S:geo}.

\begin{example}\label{Ex:exMu}
Consider the case that $\Sigma=\left\lbrace 0\right\rbrace$. Then $(\Delta, C_0^\infty(N))$ is Markov unique if and only if $d\geq 2$. See \cite[p.114]{Takeda92}.
\end{example}

We turn to essential self-adjointness. The following theorem provides a characterization in term of the $\alpha,2$-capacity of $\Sigma$. 

\begin{theorem}\label{T:esa} Let $0<\alpha\leq 2$.
The fractional Laplacian $((-\Delta)^{\alpha/2}, C_0^\infty(N))$ is essentially self-adjoint if and only if $\cpct_{\alpha,2}(\Sigma)=0$.
\end{theorem}

For the case $\alpha=2$ Theorem \ref{T:esa} is partially implied by \cite[Theorems 10.3 and 10.5]{AGMST13}, which also imply corresponding results for powers of the Laplacian of higher integer order. In \cite{HinzMasamune} we provide a version of Theorem \ref{T:esa} for the Laplacian ($\alpha=2$) on complete Riemannian manifolds, generalizing 
earlier results given in \cite[Theorem 3]{Masamune1999} and \cite[Theoreme 1]{CdV82}.

The following is a well-known guiding example for $\alpha=2$, for the case $0<\alpha <2$ see Section \ref{S:geo}.

\begin{example}\label{Ex:exesa} 
Consider the case that $\Sigma=\left\lbrace 0\right\rbrace$. Then $(\Delta, C_0^\infty(N))$ is essentially self-adjoint if and only if $d\geq 4$. See \cite[p.114]{Takeda92} and \cite[Theorem X.11, p.161]{RS80}.
\end{example}

We formulate a proof of Theorem \ref{T:esa}. Theorem \ref{T:Mu} can be obtained by similar arguments.

\begin{proof}
Suppose that $\cpct_{\alpha,2}(\Sigma)=0$. Let $(\mathcal{L}^{(\alpha)}, \dom \mathcal{L}^{(\alpha)})$ denote the closure in $L^2(\mathbb{R}^d)$ of $-(-\Delta)^{\alpha/2}$ with initial domain $C_0^\infty(N)$.
Since clearly $\dom \mathcal{L}^{(\alpha)}\subset H^\alpha(\mathbb{R}^d)$, it suffices to show the converse inclusion. Given $u\in H^\alpha(\mathbb{R}^d)$, let $(u_n)_n\subset C_0^\infty(\mathbb{R}^d)$ be a sequence approximating $u$ in $H^\alpha(\mathbb{R}^d)$. By (\ref{E:algcap}) there is a sequence $(v_k)_k\subset C_0^\infty(N)$ such that $v_k\to 0$ in $H^\alpha(\mathbb{R}^d)$ and for each $k$, $v_k$ equals one on a neighborhood of $\Sigma$. Set $w_{nk}:=(1-v_k)u_n$ to obtain functions $w_{nk}\in C_0^\infty(N)$. Let $n$ be fixed. It is easy to see that $u_n-w_{nk}=u_nv_k\to 0$ in $L^2(\mathbb{R}^d)$ as $k\to\infty$. Because the graph norm of $(-\Delta)^{\alpha/2}$ provides an equivalent norm in $H^\alpha(\mathbb{R}^d)$, it now suffices to note that 
\begin{equation}\label{E:desired}
(-\Delta)^{\alpha/2}(u_n-w_{nk})=(-\Delta)^{\alpha/2}(u_nv_k)\to 0\quad \text{in $L^2(\mathbb{R}^d)$ as $k\to\infty$}.
\end{equation}
For any $f,g\in C_0^\infty(\mathbb{R}^d)$ we can use the identity
\begin{equation}\label{E:carre}
-(-\Delta)^{\alpha/2}(fg)=2\Gamma^{(\alpha)}(f,g)-f(-\Delta)^{\alpha/2}g-g(-\Delta)^{\alpha/2}f
\end{equation}
to define the \emph{carr\'e du champ} $\Gamma^{(\alpha)}(f,g)$ of $f$ and $g$ associated with $-(-\Delta)^{\alpha/2}$, see for instance \cite[Section 1.4.2]{BGL14}. We have 
\[\left\|f(-\Delta)^{\alpha/2}g\right\|_{L^2(\mathbb{R}^d)}\leq \left\|f\right\|_{L^\infty(\mathbb{R}^d)}\left\|g\right\|_{H^\alpha(\mathbb{R}^d)}\] 
for the second summand on the right hand side, and  
\[\left\|g(-\Delta)^{\alpha/2}f\right\|_{L^2(\mathbb{R}^d)}\leq \left\|(-\Delta)^{\alpha/2}f\right\|_{L^\infty(\mathbb{R}^d)}\left\|g\right\|_{L^2(\mathbb{R}^d)}\] 
for the third. For the first summand on the right hand side of (\ref{E:carre}) we can use Cauchy-Schwarz, 
$|\Gamma^{(\alpha)}(f,g)|\leq \Gamma^{(\alpha)}(f,f)^{1/2}\Gamma^{(\alpha)}(g,g)^{1/2}$, and since $(-(-\Delta)^{\alpha/2}, C_0^\infty(\mathbb{R}^d))$ also extends to a Feller generator on $\mathbb{R}^d$ (see for instance \cite{Sa99}), we have $\Gamma^{(\alpha)}(f,f)\in L^\infty(\mathbb{R}^d)$, so that 
\[\left\|\Gamma^{(\alpha)}(f,g)\right\|_{L^2(\mathbb{R}^d)}\leq \left\|\Gamma^{(\alpha)}(f,f)\right\|_{L^\infty(\mathbb{R}^d}\left\|g\right\|_{H^\alpha(\mathbb{R}^d)}.\] 
Here we have used that $\left\|\Gamma^{(\alpha)}(g,g)\right\|_{L^1(\mathbb{R}^d)}^2$ is nothing but the energy $\left\langle (-\Delta)^{\alpha/4}g,(-\Delta)^{\alpha/4}g\right\rangle_{L^2(\mathbb{R}^d)}$ of $g$, clearly dominated by the square of the $H^\alpha(\mathbb{R}^d)$-norm of $g$. Considering (\ref{E:carre}) with $u_n$ and $v_k$ in place of $f$ and $g$ and applying the preceding estimates, we see (\ref{E:desired}). As a consequence, we see that $H^\alpha(\mathbb{R}^d)\subset \dom \mathcal{L}^{(\alpha)}$.

Conversely, suppose that $((-\Delta)^{\alpha/2}, C_0^\infty(N))$ is essentially self-adjoint in $L^2(\mathbb{R}^d)$. Then its unique self-adjoint extension must be $((-\Delta)^{\alpha/2}, H^\alpha(\mathbb{R}^d))$. Let $u\in C_0^\infty(\mathbb{R}^d)$ be a function that equals one on a neighborhood of $\Sigma$. Since $C_0^\infty(\mathbb{R}^d)\subset H^\alpha(\mathbb{R}^d)$ and by hypothesis $C_0^\infty(N)$ must be dense in $H^\alpha(\mathbb{R}^d)$, we can find a sequence $(u_n)_n$ approximating $u$ in $H^\alpha(\mathbb{R}^d)$. The functions $e_n:=u-u_n$ then are in $C_0^\infty(\mathbb{R}^d)$, equal one on a neighborhood of $\Sigma$, and converge to zero in $H^\alpha(\mathbb{R}^d)$, so that  $\cpct_{\alpha,2}(\Sigma)\leq \lim_n\left\|e_n\right\|_{H^\alpha(\mathbb{R}^d)}^2=0$.
\end{proof}

Finally, we would like to mention known removability results for $\Delta$. One says that a compact set $K\subset \mathbb{R}^d$ is \emph{removable} (or a \emph{removable singularity}) for $\Delta$ in $L^2$ if any solution $u$ of $\Delta u=0$ in $U\setminus K$ for some bounded open neighborhood $U$ of $K$ such that $u\in L^2(U\setminus K)$, can be extended to a function $\widetilde{u}\in L^2(U)$ satisfying $\Delta\widetilde{u}=0$ in $U$. See \cite[Definition 2.7.3]{AH96}. By Corollary \cite[3.3.4]{AH96} (see also \cite[Section 13.4]{Mazja} and \cite[Proposition 10.2]{AGMST13}) a compact set $K\subset \mathbb{R}^d$ is removable for $\Delta$ in $L^2$ if and only if $\cpct_{2,2}(K)=0$.

Removability results for fractional Laplacians are for instance discussed in \cite{Jayeetal16}.

\section{Riesz capacities and Hausdorff measures}\label{S:geo}

In this section we consider some geometric descriptions for the critical size of $\Sigma$. For the case of Markov uniqueness they have been discussed in many places. For the case of essential self-adjointness of integer powers of the Laplacian they were already stated in \cite{AGMST13}. 

We first give a quick review of Riesz energies and capacities. Given $s>0$ and a Radon measure $\mu$ on $\mathbb{R}^d$, let
\[I_s\mu=\int_{\mathbb{R}^d}\int_{\mathbb{R}^d}|x-y|^{-s}\mu(dy)\mu(dx)\]
denote the \emph{Riesz energy of order $s$} of $\mu$. The \emph{Riesz energy of order zero} of a Radon measure $\mu$ on $\mathbb{R}^d$ we define to be
\[I_0\mu=\int_{\mathbb{R}^d}\int_{\mathbb{R}^d} (-\ln |x-y|)_+\:\mu(dy)\mu(dx).\]
For a Borel set $E\subset\mathbb{R}^d$ we can the define the \emph{Riesz capacity of order $s\geq 0$} of $E$ by
\[\cpct_s(E)=\left[\inf\left\lbrace I_s(\mu): \text{$\mu$ Borel probability measure on $E$}\right\rbrace\right]^{-1}\]
with the agreement that $\frac{1}{\infty}:=0$. See for instance \cite[Appendix C]{Khoshnevisan02}. 

Now suppose $0<2\alpha\leq d$ and that $K\subset\mathbb{R}^d$ is compact. Then  
\begin{equation}\label{E:capiff}
\cpct_{\alpha,2}(K)>0\quad \text{ if and only if }\quad \cpct_{d-2\alpha}(K)>0.
\end{equation}
To see this note that if there exists a Borel probability measure $\mu$ on $K$ with $I_{d-2\alpha}(\mu)<+\infty$, then by (\ref{E:asympinfty}) and (\ref{E:asympzero}) respectively (\ref{E:asympzerolog}) we have $E_{\alpha}(\mu)<+\infty$, and by (\ref{E:dualdef}) therefore $\cpct_{\alpha,2}(K)>0$. Conversely, if the $\alpha,2$-capacity of $K$ is positive, we can find a nonzero Radon measure $\mu$ on $K$ with $E_\alpha(\mu)<+\infty$, so that again by (\ref{E:asympinfty}) and (\ref{E:asympzero}) respectively (\ref{E:asympzerolog}) the Borel probability measure $\frac{\mu}{\mu(K)}$ has finite Riesz energy of order $d-2\alpha$.

Consider the Dirac measure $\delta_0$ with total mass one at the origin, it is the only possible probability measure on the compact set $\left\lbrace 0\right\rbrace$. If $2\alpha\leq d$ then obviously $I_{d-2\alpha}(\delta_0)=+\infty$, so that by (\ref{E:capiff}) we have $\cpct_{\alpha,2}(\left\lbrace 0\right\rbrace)=0$. On the other hand, for $d<2\alpha$ identity (\ref{E:intcase}) implies that $U^\alpha \delta_0(x)=\gamma_{2\alpha}\ast \delta_0(x)=\gamma_{2\alpha}(x)$, $x\in\mathbb{R}^d$, so that $E_\alpha(\delta_0)=\gamma_{2\alpha}(0)<+\infty$ and therefore $\cpct_{\alpha,2}(\left\lbrace 0\right\rbrace)>0$. Similar arguments are valid with $\alpha$ in place of $2\alpha$. This produces fractional versions of Examples \ref{Ex:exMu} and \ref{Ex:exesa}.

\begin{example}
Consider the case that $0<\alpha<2$ and $\Sigma=\left\lbrace 0\right\rbrace$. Then $((-\Delta)^{\alpha/2}, C_0^\infty(N))$ is always Markov unique for $d\geq 2$. For $d=1$ it is Markov unique if $0 <\alpha\leq 1$ but not if $1<\alpha<2$. See also \cite[Section II.5, p.63]{Ber96}. So a necessary and sufficient condition for Markov uniqueness is $d\geq \alpha$.
\end{example}

\begin{example}
Consider the case that $0<\alpha<2$ and $\Sigma=\left\lbrace 0\right\rbrace$. Then $((-\Delta)^{\alpha/2}, C_0^\infty(N))$ is always essentially self-adjoint for $d\geq 4$. For $d\leq 3$ it is essentially self-adjoint if $0 <2\alpha\leq d$ but not if $d<2\alpha<4$. Therefore a necessary and sufficient condition for essential self-adjointness is $d\geq 2\alpha$.
\end{example}

As before let $\Sigma\subset\mathbb{R}^d$ be compact and of zero Lebesgue measure and write $N:=\mathbb {R}^d\setminus \Sigma$. Using theorems of Frostman-Taylor type, \cite[Appendix C, Theorems 2.2.1 and 2.3.1]{Khoshnevisan02}, see also \cite{Falconer90, Kahane85, Mattila95, MPbook}, we can give another description of the critical size of $\Sigma$, now in terms of its Hausdorff measure and dimension. 
Given $s\geq 0$, the symbol $\mathcal{H}^s$ denotes the $s$-dimensional Hausdorff measure on $\mathbb{R}^d$, \cite{Falconer90, Kahane85, Mattila95, MPbook}. By $\dim_H$ we denote the Hausdorff dimension. Again we begin with a result on Markov uniqueness.

\begin{corollary}\label{C:Mu} Let $0<\alpha\leq 2$ and suppose $\alpha\leq d$.
\begin{enumerate}
\item[(i)] If $\mathcal{H}^{d-\alpha}(\Sigma)<+\infty$ then $((-\Delta)^{\alpha/2}, C_0^\infty(N))$ is Markov unique. This is true in particular if $\alpha<d$ and $\dim_H \Sigma<d-\alpha$.
\item[(ii)] If $((-\Delta)^{\alpha/2}, C_0^\infty(N))$ is Markov unique then $\dim_H\Sigma\leq d-\alpha$.
\end{enumerate}
\end{corollary}

For the essential self-adjointness we have the following result, it partially generalizes \cite[Theorem 10.3, Corollary 10.4 and Theorem 10.5]{AGMST13}

\begin{corollary}\label{C:esa} Let $0<\alpha\leq 2$ and suppose $2\alpha\leq d$.
\begin{enumerate}
\item[(i)] If $\mathcal{H}^{d-2\alpha}(\Sigma)<+\infty$ then $((-\Delta)^{\alpha/2}, C_0^\infty(N))$ is essentially self-adjoint. This is true in particular if $2\alpha<d$ and $\dim_H \Sigma<d-2\alpha$.
\item[(ii)] If $((-\Delta)^{\alpha/2}, C_0^\infty(N))$ is essentially self-adjoint then $\dim_H\Sigma\leq d-2\alpha$.
\end{enumerate}
\end{corollary}

We provide some arguments for Corollary \ref{C:esa}, it follows from Theorem \ref{T:esa}. In a similar manner one can deduce Corollary \ref{C:Mu} from Theorem \ref{T:Mu}.
\begin{proof}
If $2\alpha<d$ and $\mathcal{H}^{d-2\alpha}(\Sigma)<+\infty$ in Corollary \ref{C:esa} (i),
then by Frostman-Taylor, \cite[Appendix C, Theorem 2.3.1]{Khoshnevisan02}, we have $\cpct_{d-2\alpha}(\Sigma)=0$, and by (\ref{E:capiff}) therefore also $\cpct_{\alpha,2}(\Sigma)=0$.
If $2\alpha=d$ and $\mathcal{H}^0(\Sigma)<+\infty$, then $\Sigma$ must be a finite set of points, note that $\mathcal{H}^0$ is the counting measure. Since capacities are subadditive, 
we have $\cpct_{d/2,2}(\Sigma)=0$ once we know a single point has zero $d/2,2$-capacity. However, a single point $p\in\mathbb{R}^d$ can only carry a point mass measure of form $c\delta_p$ where $c>0$ is a constant and $\delta_p$ the Dirac measure (with total mass one), and clearly $I_0(c\delta_p)=+\infty$, so that $\cpct_0(\left\lbrace p\right\rbrace)=0$. By (\ref{E:capiff}) this implies that $\cpct_{d/2,2}(\left\lbrace p\right\rbrace)=0$, as desired.  
Conversely, if we have $2\alpha<d$ and $\cpct_{\alpha,2}(\Sigma)=0$ Corollary \ref{C:esa} (ii), then by (\ref{E:capiff}) $\cpct_{d-2\alpha}(\Sigma)=0$, and Frostman-Taylor implies that for any $\varepsilon>0$, $\mathcal{H}^{d-2\alpha+\varepsilon}(\Sigma)=0$, showing $\dim_H\Sigma\leq d-2\alpha$. If $2\alpha=d$ and $\cpct_{d/2,2}(\Sigma)=0$, then by (\ref{E:capiff}) we have $\cpct_0(\Sigma)=0$. It is not difficult to see that this implies $\cpct_{\varepsilon}(\Sigma)=0$ for all $\varepsilon>0$, and therefore $\dim_H\Sigma=0$.
\end{proof}

\section{Additive processes and a probabilistic characterization}\label{S:Probab}

In this section we provide probabilistic characterizations of Markov uniqueness and essential self-adjointness. We use the notation $\mathbb{R}_+=[0,+\infty)$.

We are aiming only at results on hitting probabilities, so there is ambiguity what sort of stochastic process to use. Potential theory suggests to use Markov processes, and due to the group structure of $\mathbb{R}^d$ a particularly simple choice is to use certain L\'evy processes,  \cite{Ber96, Sa99}. Recall that a \emph{L\'evy process on $\mathbb{R}^d$} is a stochastic process $(X_t)_{t\in\mathbb{R}_+}$, modelled on a probability space $(\Omega,\mathcal{F},\mathbb{P})$ and taking values in $\mathbb{R}^d$ that has independent and stationary increments, is stochastically continuous, is $\mathbb{P}$-a.s. right-continuous with left limits ('c\`adl\`ag') and such that $\mathbb{P}(X_0=0)=1$. See for instance \cite[Chapter I, Section 1, Definition 1.6]{Sa99}.

Let $(B_t)_{t\in\mathbb{R}_+}$ denote a \emph{Brownian motion} on $\mathbb{R}^d$ (starting at the origin), modelled on a probability space $(\Omega,\mathcal{F},\mathbb{P})$, that is a L\'evy process on $\mathbb{R}^d$ with $\mathbb{P}$-a.s. continuous paths and such that for any $t>0$ and any Borel set $A\subset\mathbb{R}^d$,
\[\mathbb{P}(B_t\in A)=\int_A p(t,x)dx,\]
where 
\[p(t,x)=\frac{1}{(2\pi t)^{d/2}}\exp\left(-\frac{|x|^2}{2t}\right), \quad t>0,\ x\in\mathbb{R}^d.\]
Alternatively, in terms of characteristic functions, a Brownian motion is a L\'evy process on $\mathbb{R}^d$ satisfying
\[\mathbb{E}\left[\exp\left\lbrace i\left\langle \xi,B_t\right\rangle\right\rbrace\right]=\exp\left\lbrace - 2^{-1} t |\xi|^2\right\rbrace \quad t\geq 0,\ \xi\in\mathbb{R}^d.\]

More generally, given $0<\alpha\leq 2$ let $(X_t^{(\alpha)})_{t\in\mathbb{R}_+}$ denote an \emph{isotropic $\alpha$-stable L\'evy process} on $\mathbb{R}^d$, modelled on a probability space $(\Omega,\mathcal{F},\mathbb{P})$, that is a L\'evy process on $\mathbb{R}^d$ satisfying
\[\mathbb{E}\left[\exp\left\lbrace i\left\langle \xi,X^{(\alpha)}_t\right\rangle\right\rbrace\right]=\exp\left\lbrace -2^{-\alpha/2} t |\xi|^\alpha\right\rbrace \quad t\geq 0,\ \xi\in\mathbb{R}^d.\]
Obviously for $\alpha=2$ the process $(X_t^{(2)})_{t\in\mathbb{R}_+}$ is equal in law to a Brownian motion $(B_t)_{t\in\mathbb{R}_+}$. For $0<\alpha<2$ an isotropic $\alpha$-stable L\'evy process can be obtained from a Brownian motion by subordination, see \cite[Chapter 6, in particular Example 30.6]{Sa99}. For general existence results for L\'evy processes see \cite[Section I.1, Theorem 1]{Ber96} or \cite[Corollary 11.6]{Sa99}.

To prepare the discussion of related two-parameter processes below, we collect some properties. Let $0<\alpha\leq 2$. By 
\[T^{(\alpha)}_tf(x)=\mathbb{E}[f(X^{(\alpha)}_t+x)],\quad t\geq 0,\ x\in\mathbb{R}^d,\]
we can define a strongly continuous contraction semigroup $(T^{(\alpha)}_t)_{t>0}$ of Markov operators on $L^2(\mathbb{R}^d)$ (and on the space $C_\infty(\mathbb{R}^d)$ of continuous functions vanishing at infinity), they are symmetric in $L^2(\mathbb{R}^d)$. Its infinitesimal generator (in both spaces) is $-2^{-\alpha/2}(-\Delta)^{\alpha/2}$. The associated $1$-resolvent operators $R_1^{(\alpha)}=(I+2^{-\alpha/2}(-\Delta)^{\alpha/2})^{-1}$ satisfy
\[R_1^{(\alpha)}f=\int_0^\infty e^{-t}T^{(\alpha)}_tf\:dt,\]
they are bounded linear operators on $L^2(\mathbb{R}^d)$ (and on $C_\infty(\mathbb{R}^d)$). The operators $R_1^{(\alpha)}$ admit radially symmetric densities $u^{(\alpha)}$, that is
\[R_1^{(\alpha)}f(x)=\int_{\mathbb{R}^d} f(y)u_1^{(\alpha)}(x-y)dy.\]
For $0<\alpha<d$ we have 
\begin{equation}\label{E:resolvdens}
c_1|x|^{\alpha-d}\leq u_1^{(\alpha)}(x)\leq c_2|x|^{\alpha-d}
\end{equation}
whenever $|x|$ is sufficiently small, where $c_1$ and $c_2$ are two positive constants. See for instance \cite[Section 10, Lemma 3.1.1 and 3.4.1]{Khoshnevisan02}.
 
Versions of Kakutani's theorem, \cite[Section 10, Theorems 3.1.1 and 3.4.1]{Khoshnevisan02}, now allow to use Brownian motions (in case $\alpha=2$) or isotropic $\alpha$-stable L\'evy processes (in case $0<\alpha<2$) to characterize Markov uniqueness. As before, $\Sigma\subset\mathbb{R}^d$ is a compact set of zero Lebesgue measure and $N:=\mathbb{R}^d\setminus \Sigma$.

\begin{corollary}
Let $0<\alpha\leq 2$ and assume $d\geq \alpha$. The operator  $((-\Delta)^{\alpha/2}, C_0^\infty(N))$ is Markov unique if and only if for any $x\notin \Sigma$ we have
\[\mathbb{P}(\text{$\exists$ $t\in\mathbb{R}_+$ such that $X^{(\alpha)}_t+x\in \Sigma$} )=0.\]
\end{corollary}

The main aim of the present note is to point out a similar characterization for essential self-adjointness. Because their definition and structure is particularly simple, we will use two-parameter additive stable processes to describe the critical size of $\Sigma$. Let $0<\alpha\leq 2$. Given two independent isotropic $\alpha$-stable L\'evy processes $(X^{(\alpha)})_{t\in \mathbb{R}_+}$ and $(\widetilde{X}^{(\alpha)})_{t\in \mathbb{R}_+}$ on $\mathbb{R}^d$ we consider the process $(\mathcal{X}^{(\alpha)}_{\mathbf{t}})_{\mathbf{t}\in \mathbb{R}_+^2}$ defined by
\begin{equation}\label{E:additivestable}
\mathcal{X}^{(\alpha)}_{\mathbf{t}}=X^{(\alpha)}_{t_1}+\widetilde{X}^{(\alpha)}_{t_2},\quad \mathbf{t}=(t_1,t_2)\in \mathbb{R}_+^2.
\end{equation}
It is called the \emph{two-parameter additive stable process if index $\alpha$}, see \cite[Section 11.4.1]{Khoshnevisan02}. In the case $\alpha=2$ it is called the \emph{two-parameter additive Brownian motion}, we also denote it by $(\mathcal{B}_{\mathbf{t}})_{\mathbf{t}\in \mathbb{R}_+}$, where
\[\mathcal{B}_{\mathbf{t}}=B_{t_1}+\widetilde{B}_{t_2},\quad \mathbf{t}=(t_1,t_2)\in \mathbb{R}_+^2,\]
with two independent Brownian motions $(B_t)_{t\in\mathbb{R}_+}$ and $(\widetilde{B}_t)_{t\in \mathbb{R}_+}$ on $\mathbb{R}^d$. Additive stable processes or, more generally, additive L\'evy processes have been studied intensely in \cite{Khoshnevisan02, KhoshnevisanXiao04, KhoshnevisanXiao05} and follow up articles.

It seems plausible that, as two processes are added, these two-parameter processes move 'more actively' than their one-parameter versions, so they should be able to hit smaller sets with positive probability. This is indeed the case and can be used for our purpose. The next satement is a simple application of known Kakutani-type theorems for two-parameter processes, \cite[Section 11, Theorem 4.1.1]{Khoshnevisan02}.

\begin{corollary}\label{C:final}
Let $0<\alpha\leq 2$ and assume $d\geq 2\alpha$. The operator  $((-\Delta)^{\alpha/2}, C_0^\infty(N))$ is essentially self-adjoint if and only if for any $x\notin \Sigma$ we have
\[\mathbb{P}(\text{$\exists$ $\mathbf{t}\in\mathbb{R}_+^2$ such that $\mathcal{X}^{(\alpha)}_{\mathbf{t}}+x\in \Sigma$} )=0.\]
\end{corollary}

Applying Corollary \ref{C:final} with $\alpha=2$ and $d\geq 4$ we can conclude that a compact set $K\subset \mathbb{R}^d$ is removable for $\Delta$ in $L^2$ if and only if it is not hit by the additive Brownian motion with positive probability.

We collect some notions and facts related to additive stable processes and then briefly comment on the case $d=2\alpha$ in Corollary \ref{C:final} which is the only case not covered by \cite[Section 11, Proposition 4.1.1 and Theorem 4.1.1]{Khoshnevisan02}. 

One can define a two-parameter family $(\mathcal{T}_{\mathbf{t}}^{(\alpha)})_{\mathbf{t}\succ 0}$
of bounded linear operators $\mathcal{T}_{\mathbf{t}}^{(\alpha)}$ on $L^2(\mathbb{R}^d)$ (or $C_\infty(\mathbb{R}^d)$) by
\[\mathcal{T}^{(\alpha)}_{\mathbf{t}}:=T_{t_1}^{(\alpha)}T_{t_2}^{(\alpha)},\quad \mathbf{t}=(t_1,t_2)\succ 0.\] 
Here we write $(t_1,t_2)\succ (s_1,s_2)$ if $t_1>s_1$ and $t_2>s_2$. They satisfy the semigroup property $\mathcal{T}_{\mathbf{t}}^{(\alpha)}\mathcal{T}_{\mathbf{s}}^{(\alpha)}=\mathcal{T}_{\mathbf{s}
+\mathbf{t}}^{(\alpha)}$ for all $\mathbf{s},\mathbf{t}\succ 0$ and also the strong limit relation
\[\lim_{\mathbf{t}\to 0}\left\|\mathcal{T}_{\mathbf{t}}^{(\alpha)}f-f\right\|_{\sup}=0,\quad f\in C_\infty(\mathbb{R}^d),\]
and, using the density of $C_\infty(\mathbb{R}^d)\cap L^2(\mathbb{R}^d)$ in $L^2(\mathbb{R}^d)$, also for the $L^2(\mathbb{R}^d)$-norm and $f\in L^2(\mathbb{R}^d)$. By the independence of the summands in (\ref{E:additivestable}) it is not difficult to see that 
\[\mathcal{T}_{\mathbf{t}}^{(\alpha)}\mathbf{1}_A(x)=\mathbb{P}(\mathcal{X}_{\mathbf{t}}^{(\alpha)}+x\in A)\]
for all Borel sets $A\subset\mathbb{R}^d$ and starting points $x\in\mathbb{R}^d$. See for instance \cite{Hirsch95} or \cite[Sections 11.1 and 11.2]{Khoshnevisan02}. Mimicking the one-parameter case, on can introduce associated $\mathbf{1}$-resolvent operators $\mathcal{R}_{\mathbf{1}}^{(\alpha)}$ by 
\[\mathcal{R}_{\mathbf{1}}^{(\alpha)}f(x)=\int_{\mathbb{R}_+^2}e^{-(s_1+s_2)}\mathcal{T}_{\mathbf{s}}^
{(\alpha)}f(x)d\mathbf{s}.\]
Here, in accordance with the notation used above, we write $\mathbf{1}=(1,1)$. Obviously \[\mathcal{R}_{\mathbf{1}}^{(\alpha)}=R_1^{(\alpha)}R_1^{(\alpha)},\] 
and consequently the 
$\mathcal{R}_{\mathbf{1}}^{(\alpha)}$ are bounded and linear operators on $L^2(\mathbb{R}^d)$ (and $C_\infty(\mathbb{R}^d)$) and admit the densities 
\begin{equation}\label{E:densconv}
u_{\mathbf{1}}^{(\alpha)}=u_1^{(\alpha)}\ast u_1^{(\alpha)},
\end{equation}
that is
\[\mathcal{R}_{\mathbf{1}}^{(\alpha)}f(x)=\int_{\mathbb{R}^d}f(y)u_{\mathbf{1}}^{(\alpha)}(x-y)dy.\]
We provide the arguments for the special case $2\alpha=d$ in Corollary \ref{C:final}. By Giraud's lemma, \cite[Chapter 4, Proposition 4.12]{Aubin82}, together with (\ref{E:resolvdens}) and (\ref{E:densconv}), the densities $u_{\mathbf{1}}^{(d/2)}$ are continuous away from the origin and satisfy
\[u_{\mathbf{1}}^{(d/2)}(x)\leq c_3\:(-\log |x|)\]
for sufficiently small $x$, where $c_3$ is a positive constant. We also have 
\[u_{\mathbf{1}}^{(d/2)}(x)\geq c_4\:(-\log |x|)\]
for sufficiently small $x$ with a positive constant $c_4$: Given $x\in\mathbb{R}^d$ with $|x|\leq 1$, any point $y\in B(0,\frac32|x|)\setminus B(x,\frac12|x|)$ satisfies
\[\frac12|x|\leq |y|\leq \frac32|x|\quad\text{ and }\quad \frac12|x|\leq |y-x|\leq \frac52|x|.\]
Consequently, using (\ref{E:resolvdens}) and (\ref{E:densconv}), and with positive constants $c$ independent of $x$ but possibly varying from line to line,
\begin{align}
u_{\mathbf{1}}^{(d/2)}(x)&\geq c\int_{B(0,\frac32|x|)\setminus B(x,\frac12|x|)}|y-x|^{-d/2}|y|^{-d/2}dy\notag\\
&\geq c\int_{B(0,\frac12|x|)\setminus B(0,\frac12|x|^2)}|y|^{-d}dy\notag\\
&\geq c\:(-\log|x|).\notag
\end{align} 
Now an application of \cite[Section 11.3, Theorem 3.1.1]{Khoshnevisan02} yields Corollary \ref{C:final} for $2\alpha=d$.

\begin{remark}
Alternatively, one can use the $\mathbb{R}^d$-valued two-parameter Brownian sheet to characterize the essential self-adjointness of $(\Delta, C_0^\infty (N))$. A real-valued Gaussian process indexed by $\mathbb{R}_+^2$ is called a \emph{two-parameter Brownian sheet} if it has mean zero and covariance function $C(\mathbf{s},\mathbf{t})=(s_1\wedge t_1)(s_2\wedge t_2)$, $\mathbf{s},\mathbf{t}\in\mathbb{R}_+^2$.
An \emph{$\mathbb{R}^d$-valued two-parameter Brownian sheet} is a process $(\mathbb{B}_{\mathbf{t}})_{\mathbf{t}\in \mathbb{R}_+^2}$, where 
\[\mathbb{B}_{\mathbf{t}}=(\mathbb{B}_{\mathbf{t}}^1,...,\mathbb{B}_{\mathbf{t}}^d),\] 
and the components $(\mathbb{B}_{\mathbf{t}}^i)_{\mathbf{t}\in\mathbb{R}_+^2}$, $i=1,...,d$, are independent two-parameter Brownian sheets. See for instance \cite{Khoshnevisan02} or \cite{KhoshnevisanShi99}. Using the arguments of \cite{HirschSong95} one can conclude that $(\Delta, C_0^\infty (N))$ is essentially self-adjoint if and only if $\Sigma$ is polar for the two-parameter Brownian sheet, more precisely, if and only if 
\[\int_{\mathbb{R}^d}\mathbb{P}(\mathbb{B}_{\mathbf{t}}+x\in\Sigma)dt=0.\]
\end{remark}

\end{document}